\documentclass[reqno,centertags]{amsart}
\usepackage{amssymb,upref}

\renewcommand{\doteq}{:=}

\newcommand{\R}{{\mathbb{R}}}

\newcommand{\be}{\begin{equation}}
\newcommand{\ee}{\end{equation}}

\newcommand{\abs}[1]{\left|#1\right|}
\newcommand{\norm}[1]{\left\|#1\right\|}
\newcommand{\Lip}{{\mathrm{Lip}}}
\newcommand{\TV}[1]{{\rm T.V.}(#1)}

\DeclareMathOperator{\Div}{div}

\newcommand{\dx}{\, dx}

\newcommand{\lb}{\label}

\allowdisplaybreaks
\numberwithin{equation}{section}

\newtheorem{theorem}{Theorem}[section]
\newtheorem{lemma}[theorem]{Lemma}

\begin{document}

\title[Stability  of parabolic equations]{Stability of solutions of 
quasilinear\\ parabolic equations}
\author[G.\ M.\ Coclite and  H.\ Holden]{Giuseppe Maria Coclite and  
Helge Holden}
\address[Coclite]{SISSA-ISAS, via Beirut 2-4, I--34014 Trieste, Italy}
\email{coclite@sissa.it}
\address[Holden]{Department of Mathematical Sciences,
Norwegian University of
Science and Technology, Alfred Getz vei 1, NO--7491 Trondheim, Norway, 
and \newline
  Centre of Mathematics for Applications, 
  P.O.\ Box 1053, Blindern,
  N--0316 Oslo, Norway}
\email{holden@math.ntnu.no}
\urladdr{http://www.math.ntnu.no/\~{}holden/}
\date{\today}

\thanks{Partially supported by the BeMatA program of the Research
  Council of Norway and the European network  HYKE, contract
  HPRN-CT-2002-00282. The first author gratefully acknowledges the
  support and hospitality of the Centre of Mathematics for Applications.}

\begin{abstract}  We  bound the difference between solutions $u$ and $v$ of
$u_t = a\Delta u+\Div_x f+h$ and $v_t = b\Delta v+\Div_x g+k$
with initial data $\varphi$ and $ \psi$, respectively, by
$\Vert u(t,\cdot)-v(t,\cdot)\Vert_{L^p(E)}\le
A_E(t)\Vert \varphi-\psi\Vert_{L^\infty(\R^n)}^{2\rho_p}+ B(t)(\Vert a-b\Vert_{\infty}+
\Vert \nabla_x\cdot f-\nabla_x\cdot g\Vert_{\infty}+ \Vert f_u-g_u\Vert_{\infty}
+ \Vert h-k\Vert_{\infty})^{\rho_p} \abs{E}^{\eta_p}$. Here all 
functions $a$, $f$, and $h$  are
smooth and bounded, and may
depend on $u$,  $x\in\R^n$, and $t$. The functions $a$ and $h$ may in addition depend on $\nabla u$. Identical assumptions hold for the functions that determine the solutions $v$.  Furthermore, $E\subset\R^n$ is
assumed to be a bounded set, and $\rho_p$ and $\eta_p$ are fractions that
depend on $n$ and $p$. The diffusion coefficients $a$ and $b$ are assumed
to be strictly positive and the initial data are smooth.
\end{abstract}
\maketitle
\section{Introduction} \lb{s1}
We show that one can bound the difference between solutions $u$ and
$v$ of
\begin{align}
u_t& = a(t,x,u,\nabla u)\Delta u+\Div_x \big(f(t,x,u)\big)+h(t,x,u,\nabla u), 
\quad x\in \R^n,\>0<t<T,\notag \\ 
u(0,x)& = \varphi (x), \quad x\in \R^n,\label{hoved_u} \\
\intertext{and} 
v_t& = b(t,x,v,\nabla v)\Delta
v+\Div_x\big(g(t,x,v)\big)+k(t,x,v,\nabla v), 
\quad x\in \R^n,\>0<t<T,\notag \\ 
v(0,x)&= \psi (x), \quad x\in \R^n,\label{hoved_v}
\end{align}
respectively. The assumptions are that the diffusion coefficients $a$ and
$b$ are bounded from below by a strictly positive constant.  
All functions $a$, $f$, $h$, etc, as well
as the initial data $\varphi$, etc, are assumed to be smooth and
bounded.  We are interested in estimating the {\em local} $L^p$-norm
of $u(t,\cdot)-v(t,\cdot)$ over any {\em bounded} subset $E\subset\R^n$ in
terms of norm differences of the initial data as well as $a$ and $b$, etc.

In the hyperbolic case, that is, $a=b=0$, the classical result of 
Kuznetsov \cite{Kuz} and Lucier \cite{Luc} (see also \cite[Ch.~2]{HR}) reads
$$
\norm{u(t,\cdot)-v(t,\cdot)}_{L^1(\R)}\le \norm{\varphi-\psi} _{L^1(\R)}
+t\, \min\{\TV{\varphi},\TV{\psi}\}\norm{f-g}_{\Lip}
$$
in the one-dimensional case ($n=1$) where $f=f(u)$, $g=g(u)$ and
$h=k=0$.  
Here $\TV{\phi}$ denotes the total variation of the function $\phi$
and $\norm{f}_{\Lip}$  denotes the Lipschitz semi-norm. Recently, 
Bianchini and Colombo \cite{BiaCol} showed flux stability in the case 
of hyperbolic systems on the line. Indeed, they established the estimate
$$
\norm{u(t,\cdot)-v(t,\cdot)}_{L^1(\R)}\le C\, t\, \norm{Df-Dg}_{C^0(\Omega)}
$$
for solutions $u$ and $v$ of $u_t+f(u)_x=0$, $v_t+g(u)_x=0$, respectively with $u|_{t=0}=v|_{t=0}$.    
 The usual assumptions on the flux functions and the initial conditions apply, see  \cite{BiaCol}.

The dependence in $a$ of the solution $u$ of the equation
$$
u_t-\Delta a(u)=0
$$
is treated in \cite{BenCran}, assuming only that $a$ is nondecreasing, and thereby allowing degenerate diffusion.  However, no  explicit stability estimate is provided. Otto \cite{otto} studied the equation
$$
B(u)_t-\Div_x(a(\nabla u, B(u)))+h(B(u))=0
$$
with a continuous and monotone nondecreasing $B$.  Under certain assumptions he proved that
$$
\norm{B(u_1(t))-B(u_2(t))}_1\le \exp(Lt) \norm{B(u_1(0))-B(u_2(0))}_1.
$$
By extending Kru{\v z}kov's famous doubling of variables method, Bouchut and Perthame \cite{BouPer} showed that 
$$
\norm{u_1(t,\cdot)-u_2(t,\cdot)}_{L^1(\R^n)}\le \norm{u_1^0-u^0_2}_{L^1(\R^n)}
+C\, \TV{u^0_1}\sqrt{t\, \Lip(a)}
$$
when $u_j$ satisfies $u_t+\Div_x(f)=\Delta a(u_j)$ with initial data $u_j^0$, $j=1,2$. Here $a$ is assumed to be Lipschitz and nondecreasing.

Closer to the approach of this paper, Cockburn and Gripenberg \cite{CocGrip} established the estimate
$$
\norm{u_1(t,\cdot)-u_2(t,\cdot)}_{L^1(\R^n)}\le \TV{\varphi}\big(t \,\norm{f'_1-f'_2}_\infty
+4 \sqrt{tn}\norm{\sqrt{a'_1}-\sqrt{a'_2}}_\infty \big)
$$ 
for solutions $u_j$, $j=1,2$ of
$$
u_{j,t}=\Div_x (f_j)+\Delta (a_j(u_j)), \quad u_j|_{t=0}=\varphi.
$$
Allowing for explicit spatial dependence in the flux function, Evje, Karlsen, and Risebro \cite{EvKarRis,KarRis} showed stability for solutions of
$$
u_{j,t}+\Div_x(k_j(x) f_j(u))=\Delta A_j(u), \quad u_j|_{t=0}=u_j^0,
$$
in the sense that
\begin{align*}
\norm{u_1(t,\cdot)-u_2(t,\cdot)}_{L^1(\R^n)}&\le \norm{u_1^0-u^0_2}_{L^1(\R^n)} +t\, C\big( \norm{k_1-k_2}_{\infty,bv}+ \norm{f_1-f_2}_{\infty,\Lip}\big)\\
&\quad
+\sqrt{t}\, C\norm{\sqrt{A_1'}-\sqrt{A_2'}}_\infty
\end{align*}
where $\norm{\,\cdot\,}_{\infty,bv}$ and $\norm{\,\cdot\,}_{\infty,\Lip}$ is the sum of the sup-norm and the BV-norm and the sum of sup-norm and the Lipschitz norm, respectively.  Here $A_j$ is  allowed  to be degenerate.  
Karlsen and Ohlberger \cite{KarOhl} established $L^1$ contractivity of solutions of
$$
u_t+\Div_x(V(t,x) f(u))=\nabla\cdot(K(t,x)\nabla A(u)))+q(t,x,u).
$$
Recently, Chen and Karlsen \cite{CK} established the estimate
\begin{align*}
\norm{u_1(t,\cdot)-u_2(t,\cdot)}_{L^1(\R^n)}&\le \norm{u_1^0-u^0_2}_{L^1(\R^n)} 
+t\, C\norm{f'_1-f'_2}_{\infty} \\
&\quad+ \Big(t\, \norm{(\sqrt{A_1}-\sqrt{A_2})(\sqrt{A_1}-\sqrt{A_2})^\top}_\infty\Big)^{1/2}
\end{align*}
for solutions of $u_{j,t}+\Div_x f_j(u_j)=\nabla\cdot(A_j(u_j)\nabla
u_j)$ with initial data  $u_j|_{t=0}=u_j^0$.

We consider here the strictly parabolic case where the diffusion
constant is not allowed to decrease to zero. However, we allow full explicit spatial and temporal dependence in all parameters. In addition, we let the diffusion and source depend explicitly on the gradient of the unknown $u$.  All parameters, including the initial data are assumed to be smooth. Existence of regular bounded solutions is secured by classical results, see \cite{L}.  The question is to obtain explicit stability estimates. 
Our main result reads as follows. Let $u$ and $v$ denote solutions of \eqref{hoved_u} and \eqref{hoved_v}, respectively. Then
\begin{align*}
\Vert u(t,\cdot)-v(t,\cdot)\Vert_{L^p(E)}&\le
A_E(t)\Vert \varphi-\psi\Vert_{L^\infty(\R^n)}^{2\rho_p}\\
&\quad+ B(t)\Big(\Vert a-b\Vert_{L^{\infty}({\mathcal R_0})}+
\Vert \nabla_x\cdot f-\nabla_x\cdot g\Vert_{L^{\infty}({\mathcal R})}\\
&\qquad+ \Vert f_u-g_u\Vert_{L^{\infty}({\mathcal R})}
+ \Vert h-k\Vert_{L^{\infty}({\mathcal R_0})}
\Big)^{\rho_p} \abs{E}^{\eta_p},
\end{align*}
where
\begin{align*}
\rho_p&\doteq\begin{cases}
\frac{1}{2},  &\text{if $1\le p\le 2$},  \\
\frac{1}{p},  &\text{if $2< p<\infty$},
\end{cases}
\qquad
\eta_p\doteq\begin{cases}
\frac{2-p}{2p}+\frac{1}{2n},  &\text{if $1\le p\le 2$},  \\
\frac{1}{np},  &\text{if $2< p<\infty$},
\end{cases}\\
 A_E(t)&\doteq C\begin{cases}
\big(\abs{E}^{(2-p)/2p+1/2n}+\abs{E}^{1/p}\big),  &\text{if $1\le p\le 2$},  \\
(1+t^{(p-2)/p})\big(\abs{E}^{1/np}+\abs{E}^{1/p}\big),  &\text{if $2< p<\infty$},
\end{cases}\\
B(t)&\doteq C\begin{cases}
t,  &\text{if $1\le p\le 2$},  \\
(t+t^{2/p}),  &\text{if $2< p<\infty$},
\end{cases}
\end{align*}
for any bounded connected set $E\subset\R^n$ with Lipschitz
boundary. Here
$\mathcal R_0=[0,T]\times E\times[-K_1,K_1]\times[-K_2,K_2]$ and 
$\mathcal R=[0,T]\times E\times[-K_1,K_1]$.

As a particular example we note that for solutions $u$ and $v$ of
$$
u_t=a(t,x,u, \nabla u)\Delta u, \quad v_t=b(t,x,v, \nabla v)\Delta v
$$
with initial conditions $u|_{t=0}=\varphi$ and $v|_{t=0}=\psi$, we find
\begin{align*}
\Vert u(t,\cdot)-v(t,\cdot)\Vert_{L^2(E)}&\le
C(\abs{E}^{1/2n}+\abs{E}^{1/2})\norm{\varphi-\psi}_{L^\infty(\R^n)}\\
&\quad +C\,t\, \abs{E}^{1/2n}\norm{a-b}_{L^{\infty}({\mathcal R_0})}^{1/2}.
\end{align*}

Our proof is based on a homotopy argument, inspired by \cite{Bressan}.  Introducing 
\begin{align*}
u_{\theta,t}&=\big(\theta a
+(1-\theta)b\big)\Delta u_{\theta} +\Div_x\big(\theta f+(1-\theta)g\big)
 + \theta h+(1-\theta)k,  \\
u_\theta|_{t=0}& = \theta\varphi +(1-\theta)\psi,  
\end{align*}
we see that $u_0=u$ and $u_1=v$.  Thus $u_\theta$ interpolates between $u$ (for $\theta=0$) and $v$ (for $\theta=1$). The key estimate establishes that
\begin{align*}
\Vert u(t,\cdot)-v(t,\cdot)\Vert_{L^p(E)}&\equiv
{\rm dist}_{L^p(E)}\big( u(t,\cdot),v(t,\cdot)\big)\\
&\le {\rm length}_{L^p (E)}\big(u_\theta(t,\cdot)\big)\\
&=\int_0^1 \Big\Vert
\frac{\partial u_\theta}{\partial \theta}(t,\cdot)\Big\Vert_{L^p(E)}d\theta,
\end{align*}
and we establish $\theta$-independent estimates for $\norm{\partial u_\theta/\partial \theta}$.
\section{Fundamental assumptions} \lb{s2}
Fix $T>0$. Let $ u=u(t,x)$  and $v=v(t,x)$ be the bounded solution of the 
quasilinear initial value problem (see \cite{L}) 
\begin{align}
u_t& = a(t,x,u, \nabla u)\Delta u+\Div_x \big(f(t,x,u)\big)
+h(t,x,u, \nabla u),  \quad x\in \R^n,\>0<t<T,\notag  \\
u(0,x)& = \varphi (x),  \quad x\in \R^n, \label{1}\\
\intertext{and} 
v_t& = b(t,x,v,\nabla v)\Delta v+\Div_x\big(g(t,x,v)\big)
+k(t,x,v,\nabla v),  
\quad x\in \R^n,\>0<t<T,\notag  \\
v(0,x)& = \psi (x),  \quad x\in \R^n, \label{2}
\end{align}
respectively. Here
$$
f=(f_1,\dots,f_n), g=(g_1,\dots,g_n)\colon
\R\times\R^n\times\R\to \R^n,
$$
and
\begin{align*}
\Div_x\big( f(t,x,u)\big)&=\sum_{j=1}^n \big(f_j(t,x,u,\nabla u)\big)_{x_j}\\
&=\sum_{j=1}^n\Big( f_{j,x_j}+\frac{\partial f_j}{\partial u}u_{x_j}\Big)\\
&=  \nabla_x\cdot f+f_u\cdot \nabla u.
\end{align*}
Observe that $\nabla_x\cdot f$ is a scalar.  The divergence operator
$\Div_x$ always acts on the spatial variables only.  By $\nabla_{\! q} a$
(similarly for $b$, $h$, and $k$) we
denote the gradient of $a$ with respect to the final $n$ variables
(where $\nabla u$ usually sits). Our fundamental assumptions are
\begin{itemize}
\item[$({\mathcal H}_1)$] the viscous coefficients $a$ and $b$ are of class $C^3 ([0,T]\times\R^{n}\times \R \times\R^{n})$
such that
\be\label{ass1}\begin{array}{llll}
 &0 < a_* \le a (\cdot,\cdot, \cdot, \cdot) \le
a^* < \infty, \quad &\Vert a \Vert_{C^3([0,T]\times\R^{n}\times \R\times\R^{n})} \le k_1,\\
&0 < b_* \le b (\cdot,\cdot, \cdot, \cdot) \le
b^* < \infty, \quad &\Vert b \Vert_{C^3([0,T]\times\R^{n}\times \R\times\R^{n})} \le k_1,
\end{array}\ee
for some positive constants $a_*,\> a^*,\> b_*,\> b^*,\> k_1$;
\item[$({\mathcal H}_2)$] the convective terms $f$ and $g$  are of
  class $C^3([0,T]\times\R^{n}\times \R)$ and the source terms $h$ and $k$
 are of class $C^3([0,T]\times\R^{n}\times \R\times\R^{n})$ such that 
for all $i,\>j,\>l\in \{1,\dots,n\}$ and any 
$\Phi\in\{f_1,\dots,f_n,g_1,\dots,g_n,h,k\}$ the following quantities
\begin{align}
&\Big\Vert\frac{\partial \Phi}{\partial x_i}\Big\Vert_{L^\infty },
\Big\Vert\frac{\partial^2 \Phi}{\partial x_i\partial x_j}\Big\Vert_{L^\infty },
\Big\Vert\frac{\partial^3 \Phi}{\partial x_i\partial x_j\partial x_l}
\Big\Vert_{L^\infty },\label{ass2}\\
&\Big\Vert\frac{\partial \Phi}{\partial u}\Big\Vert_{L^\infty },
\Big\Vert\frac{\partial^2 \Phi}{\partial u^2}\Big\Vert_{L^\infty },
\Big\Vert\frac{\partial^2 \Phi}{\partial x_i\partial u}
\Big\Vert_{L^\infty },\Vert\nabla_{\! q} h
\Vert_{L^\infty },\Vert\nabla_{\! q} k
\Vert_{L^\infty }, \notag
\end{align}
are all bounded by  a positive constant $k_2$;
\item[$({\mathcal H}_3)$]  the initial data $\varphi$ and $\psi$  are of class
$C^2 (\R^n)$ such that
\be\label{ass3} 
\Vert\varphi\Vert_{C^2(\R^n)},\Vert\psi\Vert_{C^2(\R^n)}\le k_3,
\ee
for a positive constant $k_3$.
\end{itemize}

\begin{lemma}[{\bf ${\mathbf L^{\infty}}$-bounds on
    ${\mathbf u}$ and ${\mathbf v}$}]\label{bounds} 
     Fix $T>0$. 
By  \cite{L} 
there exist positive constants $K_1,\> K_2,\> K_3$ such that
\begin{align}
\Vert u\Vert_{L^\infty ([0,T]\times\R^n)},
\Vert v\Vert_{L^\infty ([0,T]\times\R^n)}&\le K_1,\notag\\
\Big\Vert\frac{\partial u}{\partial x_i}\Big\Vert_{L^\infty([0,T]\times\R^n) },
\Big\Vert\frac{\partial v}{\partial x_i}
\Big\Vert_{L^\infty([0,T]\times\R^n) }&\le K_2,\label{bound1}\\
\Big\Vert\frac{\partial^2 u}{\partial x_i\partial x_j}
\Big\Vert_{L^\infty ([0,T]\times\R^n)},
\Big\Vert\frac{\partial^2 v}{\partial x_i\partial x_j}
\Big\Vert_{L^\infty([0,T]\times\R^n) }&\le K_3,\notag
\end{align}
for  all $i,\>j\in \{1,...,n\}$ where $K_1,\> K_2,\> K_3$ depend only on $T,\>  n,\> a_*,\>a^*,\>b_*,\>b^*, \> k_1,$
$k_2$, and $k_3$.
\end{lemma}

\section[]{The homotopy argument}\label{section2}

Our approach is based on the following homotopy argument. 
Let $0\le \theta\le 1$. The function $u_\theta$ interpolates between 
the functions $u$ and $v$. More precisely, denote by  $u_\theta$  the 
solution of the  quasilinear initial value problem
\begin{align}
u_{\theta,t}&=\big(\theta a(t,x,u_\theta,\nabla u_\theta)
+(1-\theta)b(t,x,u_\theta,\nabla u_\theta)\big)\Delta u_{\theta}\notag\\
&\quad +\Div_x\big(\theta f(t,x,u_\theta)+(1-\theta)g(t,x,u_\theta)\big)\notag\\
&\quad + \theta h(t,x,u_\theta,\nabla u_\theta)+(1-\theta)k(t,x,u_\theta,\nabla u_\theta),
\quad x\in \R^n,\>0<t<T,\notag  \\
u_\theta(0,x)& = \theta\varphi (x)+(1-\theta)\psi (x),  
\quad x\in \R^n.\label{teta}
\end{align}
Clearly
$$
u_0= v,\qquad u_1= u.
$$
Indeed
$$
\theta\longmapsto u_\theta(t,\cdot)
$$
is a curve joining $v(t,\cdot)$ and $u(t,\cdot)$, and
\be\label{curv1}
\Vert u(t,\cdot)-v(t,\cdot)\Vert_{L^p(E)}\equiv
{\rm dist}_{L^p(E)}\big( u(t,\cdot),v(t,\cdot)\big)\le {\rm length}_{L^p(E)}\big(u_\theta(t,\cdot)\big),
\ee
for each $0\le t\le T$, $E\subset \R^n$ measurable set  and $1\le p\le \infty$.

\begin{lemma}[{\bf ${\mathbf L^{\infty}}$-bounds on ${\mathbf u_\theta}$}]
\label{bounds'} By
\cite[Theorem V 8.1]{L}, there exist positive constants $K_1,\> K_2,\> K_3$  depending only on 
$T$,  $n,\> a_*,\>a^*,\>b_*,\>b^*,\> k_1,$
$k_2$ and $k_3$ such that
\begin{align}
\Vert u_\theta\Vert_{L^\infty ([0,T]\times\R^n)}&\le K_4,\notag \\
\Big\Vert\frac{\partial u_{\theta}}{\partial x_i}\Big\Vert_{L^\infty ([0,T]\times\R^n)}&\le K_5,\label{bound1'} \\
\Big\Vert\frac{\partial^2 u_{\theta}}{\partial x_i\partial x_j}\Big\Vert_{L^\infty ([0,T]\times\R^n)}&\le K_6,\notag
\end{align}
for each $0\le \theta\le 1$ and $i,\>j\in \{1,...,n\}$.
\end{lemma}
\begin{lemma}[{\bf Smoothness of ${\mathbf \theta \mapsto
      {u}_\theta}$}] \label{smooth} 
 Assume $({\mathcal H}_1)$, $({\mathcal H}_2)$, and $({\mathcal H}_3)$. The curve
$$
\theta\in [0,1]\longmapsto u_\theta(t,\cdot)\in C^2(\R^n)
$$
is of class $C^1$. In particular, we infer
\be\label{curv2}
{\rm length}_{L^p (E)}\big(u_\theta(t,\cdot)\big)=\int_0^1 \Big\Vert
\frac{\partial u_\theta}{\partial \theta}(t,\cdot)\Big\Vert_{L^p(E)}d\theta,
\ee
for each $0\le t\le T$ and $E\subset \R^n$ measurable set.
\end{lemma}
\begin{proof} Consider the map
\begin{align*}
F&\colon{\mathcal D}\longrightarrow C^{\infty}\big(]0,T[\times\R^n\big)\cap
C^2\big([0,T]\times\R^n\big),\\
F(\theta,\omega)(t,x)&\doteq
\frac{\partial \omega}{\partial t}(t,x)-
\Big(\theta a\big(t,x,\omega(t,x),\nabla \omega(t,x)\big)\\
&\qquad\qquad\qquad\quad +(1-\theta)b\big(t,x,\omega(t,x),\nabla \omega(t,x)\big)\Big)
\Delta \omega(t,x)\\
&\quad-\Div_x\big(\theta f(t,x,\omega(t,x))+(1-\theta)g(t,x,\omega(t,x))\big)\\
&\quad-\big(\theta h(t,x,\omega(t,x),\nabla \omega(t,x))+(1-\theta)k(t,x,\omega(t,x),\nabla \omega(t,x))\big),
\end{align*}
where
$$
{\mathcal D}\doteq\Big\{(\theta,\omega)\in [0,1]\times C^{\infty}\big(]0,T[\times\R^n\big)\cap
  C^2\big([0,T]\times\R^n\big)\mid  \omega(0,\cdot)= \theta\varphi+(1-\theta)\psi \Big\}.
$$
 From the definition of $u_\theta$,
\be\label{curv3}
F(\theta,u_\theta)\equiv 0,\qquad 0\le \theta \le 1.
\ee
Observe that $F$ is of class $C^1$ and
\begin{align*}
\frac{\partial F}{\partial \theta}(\theta,\omega)&=
\big(b(t,x,\omega,\nabla \omega)-a(t,x,\omega,\nabla \omega)\big)\Delta \omega\\
&\quad+\Div_x\big(g(t,x,\omega)-f(t,x,\omega)\big)
+  k(t,x,\omega,\nabla \omega)-h(t,x,\omega,\nabla \omega).
\end{align*}
To compute 
$$
\frac{\partial F}{\partial \omega}(\theta,\omega)\big[(\theta',z)\big]=
\left. \frac{\partial F}{\partial \varepsilon}(\theta,\omega+\varepsilon z)\right\vert_{\varepsilon=0}
$$
we find
\begin{align*}
F(\theta,\omega+\varepsilon z)&=
\frac{\partial \omega}{\partial t}+
\varepsilon\frac{\partial z}{\partial t}-
\big(\theta a\big(t,x,\omega+\varepsilon z,\nabla\omega+\varepsilon\nabla z)\\
&\qquad\qquad+(1-\theta)b(t,x,\omega+\varepsilon z,\nabla\omega+\varepsilon\nabla z)\big)
\left(\Delta \omega+\varepsilon\Delta z\right)\\
&\quad -\Div_x\big(\theta f(t,x,\omega+\varepsilon z)+
(1-\theta)g(t,x,\omega+\varepsilon z)\big)\\
&\quad -\big(\theta h(t,x,\omega+\varepsilon
z,\nabla\omega+\varepsilon\nabla
z)+(1-\theta)k(t,x,\omega+\varepsilon z, 
\nabla\omega+\varepsilon\nabla z)\big),\\
\frac{\partial F}{\partial \varepsilon}(\theta,\omega+\varepsilon z)&=
\frac{\partial z}{\partial t}-\big(\theta a(t,x,\omega+\varepsilon
z,\nabla\omega+\varepsilon\nabla z)\\
&\qquad +(1-\theta)b(t,x,\omega+
\varepsilon z,\nabla\omega+\varepsilon\nabla z)\big)
\Delta z\\
&\quad -
\Big(\theta\frac{\partial a}{\partial \omega} (t,x,\omega+\varepsilon
z,\nabla\omega+\varepsilon\nabla z)\\
&\qquad +
(1-\theta)\frac{\partial b}{\partial \omega} (t,x,\omega+\varepsilon z,\nabla\omega+\varepsilon\nabla z)\Big)
z \left( \Delta \omega+\varepsilon\Delta z\right)\\
&\quad -\Big(\theta \nabla_{\! q} a(t,x,\omega+\varepsilon
  z,\nabla\omega+\varepsilon\nabla z)\\
&\qquad +(1-\theta)\nabla_{\! q} b(t,x,\omega+\varepsilon z,
\nabla\omega+\varepsilon\nabla z)\Big)\cdot\nabla z\,
 \left( \Delta \omega+\varepsilon\Delta z\right)\\
&\quad -
\Div_x\left(
\Big(\theta\frac{\partial f}{\partial \omega}(t,x,\omega+\varepsilon z)+
(1-\theta)\frac{\partial g}{\partial \omega}(t,x,\omega+\varepsilon z)\Big)z\right)\\
&\quad -
\Big(\theta\frac{\partial h}{\partial \omega}(t,x,\omega+\varepsilon
z,\nabla\omega+\varepsilon\nabla z)\\
&\qquad+(1-\theta)\frac{\partial k}{\partial \omega}
(t,x,\omega+\varepsilon z,\nabla\omega+\varepsilon\nabla z)\Big)z\\
&\quad -
\Big(\theta\nabla_{\! q} h(t,x,\omega+\varepsilon
z,\nabla\omega+\varepsilon\nabla z)\\
&\qquad+(1-\theta)\nabla_{\! q} k(t,x,\omega+\varepsilon z,
\nabla\omega+\varepsilon\nabla z)\Big)\cdot\nabla z.
\end{align*}
Thus
\begin{align*}
\frac{\partial F}{\partial \omega}(\theta,\omega)\big[(\theta',z)\big]&=
\frac{\partial z}{\partial t}-
\big(\theta a(t,x,\omega,\nabla\omega)+(1-\theta)b(t,x,\omega,\nabla\omega)\big)\Delta z \\
&\quad-\left(\theta\frac{\partial a}{\partial \omega}(t,x,\omega,\nabla\omega)+(1-\theta)
\frac{\partial b}{\partial \omega}(t,x,\omega,\nabla\omega)\right)z\Delta\omega\\
&\quad-\left(\theta\nabla_{\! q} a(t,x,\omega,\nabla\omega)+(1-\theta)
\nabla_{\! q} b(t,x,\omega,\nabla\omega)\right)\cdot \nabla z\,\Delta\omega\\
&\quad
-\Div_x\left(\Big(\theta \frac{\partial f}{\partial \omega}(t,x,\omega)
+(1-\theta)\frac{\partial g}{\partial \omega}(t,x,\omega)\Big)z\right)\\
&\quad
- \Big(\theta \frac{\partial h}{\partial
  \omega}(t,x,\omega,\nabla\omega)
+(1-\theta)
\frac{\partial k}{\partial \omega}(t,x,\omega,\nabla\omega)\Big)z\\
&\quad
- \Big(\theta \nabla_{\! q} h(t,x,\omega,\nabla\omega)+(1-\theta)
\nabla_{\! q} k(t,x,\omega,\nabla\omega)\Big)\cdot\nabla z,\\
&\qquad\qquad (\theta,\omega),(\theta',z)\in{\mathcal D}.
\end{align*}
Observe that $(\theta',z)\in {\mathcal D}$ satisfies the equation
$$\frac{\partial F}{\partial \omega}(\theta,\omega)\big[(\theta',z)\big]=\zeta$$
if and only if $z$ is  solution of the  linear  initial value problem
\begin{align*}
z_{t}& =\big(\theta a(t,x,\omega,\nabla\omega)+(1-\theta)b(t,x,\omega,\nabla\omega)\big) \Delta z\\
&\quad+\big(\theta  a_\omega(t,x,\omega,\nabla\omega)+(1-\theta)
 b_\omega(t,x,\omega,\nabla\omega)\big) \Delta\omega z\\
 &\quad+\big(\theta  \nabla_{\! q} a(t,x,\omega,\nabla\omega)+(1-\theta)
 \nabla_{\! q} b(t,x,\omega,\nabla\omega)\big)\cdot\nabla z\, \Delta\omega \\
 &\quad+\Div_x \Big(\big(\theta f_\omega(t,x,\omega)+
 (1-\theta)g_\omega(t,x,\omega)\big)z\Big)\\
 &\quad+\big(\theta h_\omega(t,x,\omega,\nabla\omega)+
 (1-\theta)k_\omega(t,x,\omega,\nabla\omega)\big)z\\
&\quad+\big(\theta \nabla_{\! q} h(t,x,\omega,\nabla\omega)+
 (1-\theta)\nabla_{\! q} k(t,x,\omega,\nabla\omega)\big)\cdot\nabla z+\zeta(t,x),\\
&\qquad \qquad x\in \R^n,\>0<t<T,  \\
z(0,x)& = \theta'\varphi (x)+(1-\theta')\psi(x),  \quad x\in \R^n.
\end{align*}
Since this problem is well-posed (see \cite[Theorem IV 5.1]{L}),
$\frac{\partial F}{\partial \omega}(\theta,\omega)$
 is invertible.
 By the implicit function theorem, the curve
$\theta\longmapsto u_\theta$ is of class $C^1$ and clearly \eqref{curv2} holds.
This concludes the proof.
\end{proof}
Differentiating  equation  \eqref{teta} with respect to $\theta$, we have
\begin{align}
\frac{\partial^2 u_\theta}{\partial t \partial \theta}&=
\big(\theta a(t,x,u_\theta,\nabla u_\theta)
+(1-\theta)b(t,x,u_\theta,\nabla u_\theta)\big)\Delta\left(
\frac{\partial u_\theta}{ \partial  \theta}\right)\notag\\
&\quad+\left(\theta \frac{\partial a}{\partial u}(t,x,u_\theta,\nabla u_\theta)
+(1-\theta)\frac{\partial b}{\partial u} (t,x,u_\theta,\nabla u_\theta)\right)
\Delta u_\theta\,
\frac{\partial u_\theta}{\partial  \theta}\notag\\
&\quad+\left(\theta \nabla_{\! q} a(t,x,u_\theta,\nabla u_\theta)
+(1-\theta)\nabla_{\! q} b(t,x,u_\theta,\nabla u_\theta)\right)\cdot
\nabla\left(\frac{\partial u_\theta}{\partial  \theta}\right)
\, \Delta u_\theta\,\notag \\
&\quad+\big(a(t,x,u_\theta,\nabla u_\theta)-b(t,x,u_\theta,\nabla u_\theta)\big)\Delta u_\theta \label{curv4}\\
&\quad+\left(\frac{\partial f}{\partial u}(t,x,u_\theta)-
\frac{\partial g}{\partial u}(t,x,u_\theta)\right)\nabla u_\theta
+\nabla_x\cdot f(t,x,u_\theta)-\nabla_x\cdot g(t,x,u_\theta)\notag\\
&\quad+\left(\theta\frac{\partial f}{\partial u}(t,x,u_\theta)+(1-\theta)
\frac{\partial g}{\partial u}(t,x,u_\theta)\right)\cdot
\nabla\left(\frac{\partial u_\theta}{\partial  \theta}\right)\notag\\
&\quad+\left(\theta\nabla_x\cdot\frac{\partial f}{\partial u}(t,x,u_\theta) 
+(1-\theta)\nabla_x\cdot\frac{\partial g}{\partial u}(t,x,u_\theta)\right)
\frac{\partial u_\theta}{\partial  \theta}\notag\\
&\quad+\left(\theta\frac{\partial^2 f}{\partial u^2}(t,x,u_\theta)+(1-\theta)
\frac{\partial^2 g}{\partial u^2}(t,x,u_\theta)\right)\cdot
\nabla u_\theta\,\frac{\partial u_\theta}{\partial  \theta}\notag\\
&\quad+h(t,x,u_\theta,\nabla u_\theta) -k(t,x,u_\theta,\nabla
u_\theta)\notag\\
&\quad+
\left(\theta\frac{\partial h}{\partial u}(t,x,u_\theta,\nabla u_\theta)  
+(1-\theta)
\frac{\partial k}{\partial u}(t,x,u_\theta,\nabla u_\theta)\right)
\frac{\partial u_\theta}{\partial  \theta}\notag \\
&\quad+
\left(\theta\nabla_{\! q} h(t,x,u_\theta,\nabla u_\theta)  
+(1-\theta)
\nabla_{\! q} k(t,x,u_\theta,\nabla u_\theta)\right)
\cdot\nabla\left(\frac{\partial u_\theta}{\partial\theta}\right). \notag
\end{align}
Denoting
\begin{align*}
z_\theta(t,x)&\doteq\frac{\partial u_\theta}{\partial  \theta},\\
\alpha(\theta,t,x)&\doteq \theta a+(1-\theta)b,\\
\beta(\theta,t,x)&\doteq \big(\theta\nabla_{\! q} a+(1-\theta)\nabla_{\! q} b
\big)\Delta u_\theta+\theta\frac{\partial f}{\partial u}+
(1-\theta)\frac{\partial g}{\partial u}\\
&\quad + \big(\theta\nabla_{\! q} h+(1-\theta)\nabla_{\! q} k \big),\\
\gamma(\theta,t,x)&\doteq\left(\theta \frac{\partial a}{\partial u}
+(1-\theta)\frac{\partial b}{\partial u} \right)
\Delta u_\theta\\
&\quad+\theta\nabla_x\cdot\frac{\partial f}{\partial u}
+(1-\theta)\nabla_x\cdot\frac{\partial g}{\partial u}
+\theta\frac{\partial h}{\partial u}+(1-\theta)
\frac{\partial k}{\partial u}\\
&\quad+
\left(\theta\frac{\partial^2 f}{\partial u^2}+(1-\theta)
\frac{\partial^2 g}{\partial u^2}\right)
\cdot\nabla u_\theta,\\
\sigma(\theta,t,x)&\doteq
\big(a-b\big)\Delta u_\theta+
\left(\frac{\partial f}{\partial u}
-\frac{\partial g}{\partial u}\right)
\cdot\nabla u_\theta\\
&\quad+\nabla_x\cdot f-\nabla_x\cdot g +h-k,
\end{align*}
for each $0\le \theta\le 1,\> t\ge 0,\>x\in \R^n,$
there results
\begin{multline}\label{teta'}
\frac{\partial z_\theta}{\partial t}=
\alpha(\theta,t,x)\Delta z_\theta+\beta(\theta,t,x)\cdot\nabla z_\theta+
\gamma(\theta,t,x)z_\theta+\sigma(\theta,t,x),\\
 0\le \theta\le 1,\>\> 0<t<T,\>\> x\in \R^n.
\end{multline}
Moreover, observe that
\be\label{teta'0}
z_\theta(0,x)=\varphi(x)-\psi(x),\qquad 0\le \theta\le 1,\>\> x\in \R^n.
\ee
\begin{lemma}[{\bf ${\mathbf L^{\infty}}$-bounds on ${\mathbf \alpha}$,
${\mathbf \beta}$, ${\mathbf \gamma}$}]\label{bounds''}
 From the definition of $\alpha$,
\eqref{ass1} and \eqref{bound1'}, we have
\be\label{b1}
0<\alpha_*\le \alpha(\cdot,\cdot,\cdot)\le \alpha^*,\qquad
\Vert \nabla\alpha\Vert_{L^\infty}\le k_1(1+ K_5+nK_6)
\ee
where
$$
\alpha_*\doteq {\rm min}\{a_*,b_*\},\qquad \alpha^*\doteq{\rm max}\{a^*,b^*\}.
$$
Moreover, from the definition of $\beta$ and \eqref{ass2}, we infer
\be\label{b2}
\Vert \beta\Vert_{L^{\infty}}= \sup_{j=1,\dots,n}\Vert
\beta_j\Vert_{L^{\infty}}\le K_7,
\ee
where
$$
K_7\doteq n k_1 K_6+2k_2.
$$
Finally, from the definition of $\gamma$, \eqref{ass1}, \eqref{ass2} and \eqref{bound1'}, we find
\be\label{b3}
\Vert \gamma\Vert_{L^{\infty}}\le K_8,
\ee
where
$$
K_8\doteq nk_1 K_6+(n+1+nK_5)k_2.
$$
\end{lemma}
\begin{lemma}[{\bf ${\mathbf L^{\infty}}$-bounds on ${\mathbf
      z_\theta}$}]\label{teta1} 
Assume
 $({\mathcal H}_1)$, $({\mathcal H}_2)$, and $({\mathcal H}_3)$.
 There exists a positive constant $C_1$ depending only on $T$,
  $n,\> a_*,\>a^*,\>b_*,\>b^*,\> k_1,\> k_2$, and $k_3$
such that
\be\label{teta'b}
\Vert z_\theta (t,\cdot)\Vert_{L^{\infty}(\R^n)}\le
C_1 t+\Vert \varphi-\psi \Vert_{L^{\infty}(\R^n)},
\ee
for each $0\le t\le T$ and $0\le \theta\le 1$.
\end{lemma}
\begin{proof} To simplify the notation we let $w$ denote the solution of \eqref{teta'}, that is,
\begin{equation}\label{st_teta'}
w_t=
\alpha\Delta w+\beta\cdot\nabla w+
\gamma w+\sigma,\quad w|_{t=0}=w_0.
\end{equation}
Linearity implies that
$$
w=w_1+w_2
$$
where $w_1$ and $w_2$ solve
\begin{align*}
w_{1,t}&=
\alpha\Delta w_1+\beta\cdot\nabla w_1+
\gamma w_1,\quad w_1|_{t=0}=w_0, \\
w_{2,t}&=
\alpha\Delta w_2+\beta\cdot\nabla w_2+
\gamma w_2+\sigma,\quad w_2|_{t=0}=0,
\end{align*}
respectively. We infer from \cite[p. 389]{L} that
\begin{align*}
w_1(t,x)&=\int_{\R^n} G(t,0,x,\xi)w_0(\xi)\, d\xi, \\
w_2(t,x)&=\int_0^t \int_{\R^n} G(t,\tau,x,\xi)\sigma(\tau,\xi)\, d\xi d\tau,
\end{align*}
where $G$ is the Green's function. For $t\in[0,T]$ for some fixed $T$ positive we find
\begin{equation*}
\abs{w_2(t,x)}\le C\, t \norm{\sigma}_\infty.
\end{equation*}
Introduce $z=w_1-w_0$ which satisfies the equation for $w_2$ with $\sigma=\alpha\Delta w_0+\beta\Div_x w_0+
\gamma w_0$. Thus
\begin{equation*}
\abs{w_1(t,x)}\le \abs{z(t,x)}+\abs{w_0(x)}\le \norm{w_0}_\infty
+C\, t \norm{\alpha\Delta w_0+\beta\cdot\nabla w_0+
\gamma w_0}_\infty.
\end{equation*}
\end{proof}
Observe that in the previous lemma, the smoothness of the initial
condition enters in a crucial way.  With less regularity we get the
familiar $O(t^{1/2})$ behavior near $t=0$ (see, e.g., \cite[Sec.~4.4]{HR}).

\section[]{Stability of quasilinear parabolic equations}\label{section2A}
We begin with the following lemma.
\begin{lemma}[{\bf Poincar\'e-type inequality}]\label{lemmaA1}
There exists a positive constant $\Lambda_0$, depending only on $n$, such that
\be\label{A1}
\int_{B} |f|^2\dx\le
\Lambda_0 \abs{B}^{2/n}
\int_{B} |\nabla f|^2\dx+
\Lambda_0 \abs{B}^{1/n}\int_{\partial  B} |f|^2\dx,
\ee
for each $f\in C^2(\R^n)$ and $B\subset \R^n$ bounded connected set
with Lipschitz boundary.
In the case $n=1$ we mean
$$
\int_{\partial  B} |f|^2\dx=\vert f(x_0)\vert^2,
$$
for some $x_0\in B$.
\end{lemma}
The proof of this lemma is more or less classical (see \cite[Theorem A.9]{S}
and \cite[Lemma A.2]{C}) and the dependence of the coefficients on the 
measure of the domain is
consequence of a standard rescaling argument.

Now we  prove the key estimate in the $L^2$-norm for the map $z_\theta$.
\begin{lemma}[{\bf Case $\mathbf {p=2}$: Energy estimate}]\label{En1}
  Assume  $({\mathcal H}_1)$, $({\mathcal H}_2)$, and $({\mathcal H}_3)$.
 Then
there exists a positive constant $C_2$  depending only on 
$T, \>n,  \> a_*,\>a^*,\>b_*,\>b^*,\> k_1,\> k_2$, and $k_3$ such that
\begin{align}
\Vert z_\theta(t,\cdot)\Vert_{L^2(E)}&\le
C_2(\abs{E}^{1/2n}+\abs{E}^{1/2})\Vert \varphi-\psi\Vert_{L^\infty(\R^n)}\notag\\
&\quad+ C_2\, t\, \abs{E}^{1/2n}\Big(\Vert a-b\Vert_{L^{\infty}({\mathcal R_0})}+
\Vert \nabla_x\cdot f-\nabla_x\cdot g\Vert_{L^{\infty}({\mathcal R})}
\label{En21}\\
&\qquad+ \Vert f_u-g_u\Vert_{L^{\infty}({\mathcal R})}
+ \Vert h-k\Vert_{L^{\infty}({\mathcal R_0})}\Big)^{1/2},\notag
\end{align}
for each $0\le t\le T$, $0\le \theta\le 1$ and $E\subset \R^n$ bounded 
connected set with Lipschitz boundary. Here  
$\mathcal R_0=[0,T]\times E\times[-K_1,K_1]\times[-K_2,K_2]$ and 
$\mathcal R=[0,T]\times E\times[-K_1,K_1]$.
\end{lemma}
\begin{proof} Let $B\subset  \R^n$ be a ball and $0<t<T$. Then
by \eqref{teta'} we find
\begin{align}\label{E1}
\frac{d}{dt}\int_{B}\frac{1}{2}z_\theta^2(t,x)\dx&=
\int_{B}z_\theta z_{\theta,t}\dx\\
&=\int_{B}\alpha z_\theta \Delta z_{\theta}\dx+
\int_{B}z_\theta\beta\cdot\nabla z_{\theta}\dx+
\int_{B}\gamma z_\theta^2 \dx+
\int_{B}\sigma z_\theta \dx.\notag
\end{align}
Observe that, by   \eqref{b3},
\be\label{E2}
\int_{B}\gamma z_\theta^2 \dx\le K_8
\int_{B} z_\theta^2 \dx,
\ee
and, by \eqref{b2},
\begin{align}
\int_{B}z_\theta \beta \cdot\nabla z_{\theta}\dx&\le
\frac{1}{\alpha_*}\int_{B}\abs{\beta}^2 z_\theta^2 \dx+
\frac{\alpha_*}{4}\int_{B}\abs{\nabla z_{\theta}}^2 \dx\label{E11}\\
&\le\frac{K_7^2}{\alpha_*}\int_{B} z_\theta^2 \dx+
\frac{\alpha_*}{4}\int_{B}\abs{\nabla z_{\theta}}^2\dx. \notag
\end{align}
By Lemma \ref{teta1} and \eqref{ass3},
\begin{align}
\int_{B}\sigma z_\theta \dx&\le
\int_{B}|\sigma|\,| z_\theta| \dx=
\int_{B}\sqrt{|\sigma|}\big(\sqrt{|\sigma|}\, |z_\theta|\big) \dx\notag\\
&\le\frac{1}{2}\int_{B}|\sigma|\dx+
\frac{1}{2}\int_{B}|\sigma|z_\theta^2 \dx\label{E3}\\
&\le\frac{1}{2}\abs{B}\,\Vert \sigma\Vert_{L^{\infty}({\mathcal R})}+
C_1^2\abs{B}t^2\Vert \sigma\Vert_{L^{\infty}({\mathcal R})}\notag\\
&\qquad\qquad+
\Vert \varphi-\psi\Vert_{L^\infty(\R^n)}^2\abs{B}\,\Vert \sigma\Vert_{L^{\infty}({\mathcal R})}\notag\\
&\le\frac{K_9}{2}\abs{B}\,\Vert \sigma\Vert_{L^{\infty}({\mathcal R})}+
C_1^2\abs{B}t^2\Vert \sigma\Vert_{L^{\infty}({\mathcal R})},\notag
\end{align}
where
$$
K_9\doteq 1+8k_3^2.
$$
Moreover,  by the divergence theorem  we have
\begin{align}
\int_{B}
\alpha z_\theta \Delta z_{\theta}\dx
&=\int_{\partial B}\alpha z_\theta (\nabla z_{\theta}\cdot\nu)\dx-
\int_{B}(\nabla \alpha\cdot \nabla z_\theta) z_{\theta}\dx-
\int_{B}\alpha  \abs{\nabla z_{\theta}}^2\dx\notag\\
&\le\int_{\partial B}\alpha z_\theta (\nabla z_{\theta}\cdot\nu)\dx+
\frac{1}{2\alpha_*}\int_{B}\abs{\nabla \alpha}^2\, z_\theta^2\dx\label{E4}\\
&\qquad
+\frac{\alpha_*}{2}\int_{B}\abs{\nabla z_{\theta}}^2\dx-
\alpha_*\int_{B}\abs{\nabla z_{\theta}}^2\dx\notag\\
&=\int_{\partial B}\alpha z_\theta (\nabla z_{\theta}\cdot\nu)\dx+
\frac{1}{2\alpha_*}\int_{B}\abs{\nabla\alpha}^2\, z_\theta^2\dx
-\frac{\alpha_*}{2}\int_{B}\abs{\nabla z_{\theta}}^2\dx\notag \\
&\le\int_{\partial B}\alpha z_\theta (\nabla z_{\theta}\cdot\nu)\dx+
\frac{\Vert\nabla\alpha\Vert^2_{L^\infty}}{2\alpha_*}\int_{B} z_\theta^2\dx
-\frac{\alpha_*}{2}\int_{B}\abs{\nabla z_{\theta}}^2\dx,\notag
\end{align}
where $\nu$ is the external normal to $\partial B$
and in the case $n=1$, $\partial B=\{x_1, x_2\},\> x_1<x_2$, we mean
$$
\int_{\partial B}\alpha z_\theta (\nabla z_{\theta}\cdot\nu)\dx=
\alpha(\theta, t, x_2) z_\theta(t, x_2) z_{\theta,x}(t, x_2)-
\alpha(\theta, t, x_1) z_\theta(t, x_1) z_{\theta,x}(t, x_1).
$$
Substituting \eqref{E2}, \eqref{E11}, \eqref{E3}, \eqref{E4} in \eqref{E1}
we obtain
\begin{align}
\frac{d}{dt}\int_{B}\frac{1}{2}z_\theta^2(t,x)\dx&\le
-\frac{\alpha_*}{4}\int_{B}\abs{\nabla z_{\theta}}^2\dx+
\int_{\partial B}\alpha z_\theta (\nabla z_{\theta}\cdot\nu)\dx\notag\\
&\quad+
\bigg(K_8+\frac{K_7^2}{\alpha_*}+
\frac{\Vert\nabla\alpha\Vert^2_{L^\infty}}{2\alpha_*}\bigg)
\int_{B} z_\theta^2\dx\label{E111}\\
&\quad+\frac{K_9}{2}\abs{B}\,\Vert \sigma\Vert_{L^{\infty}({\mathcal R})}+
C_1^2\abs{B}t^2\Vert \sigma\Vert_{L^{\infty}({\mathcal R})}.\notag
\end{align}
By Lemma \ref{lemmaA1} and the assumptions on $B$,
\be\label{E411}
-\int_{B}|\nabla z_{\theta}|^2\dx\le
-\frac{1}{\Lambda_0 \abs{B}^{2/n}}\int_{B}z_{\theta}^2\dx+
\frac{1}{\abs{B}^{1/n}}
\int_{\partial  B} z_{\theta}^2\dx,
\ee
so by Lemma \ref{teta1},  \eqref{E111} and \eqref{E411},
\begin{align}
\frac{d}{dt}\int_{B}\frac{1}{2}z_\theta^2(t,x)\dx 
&\le\bigg(K_8+\frac{K_7^2}{\alpha_*}+
\frac{\norm{\nabla\alpha}^2_{L^\infty}}{2\alpha_*}-
\frac{\alpha_*}{4 \Lambda_0 \abs{B}^{2/n}}\bigg)
\int_{B}z_{\theta}^2\dx\notag\\
&\quad+\frac{K_9\abs{B}}{2}\Vert \sigma\Vert_{L^{\infty}({\mathcal R})}+
 C_1^2\abs{B}t^2\Vert \sigma\Vert_{L^{\infty}({\mathcal R})}\notag\\
&\quad+
\int_{\partial B}\alpha z_\theta (\nabla z_{\theta}\cdot\nu)\dx
+\frac{\alpha_*}{4 \abs{B}^{1/n}}
\int_{\partial  B} z_{\theta}^2\dx\label{E5}\\
&\le\bigg(K_8+\frac{K_7^2}{\alpha_*}+
\frac{\Vert\nabla\alpha\Vert^2_{L^\infty}}{2\alpha_*}-
\frac{\alpha_*}{4\Lambda_0\abs{B}^{2/n}}\bigg)
\int_{B}z_{\theta}^2\dx\notag \\
&\quad+\frac{K_9\abs{B}}{2}\Vert \sigma\Vert_{L^{\infty}({\mathcal R})}
+C_1^2\abs{B}t^2\Vert \sigma\Vert_{L^{\infty}({\mathcal R})}\notag \\
&\quad+\int_{\partial B}\alpha z_\theta (\nabla z_{\theta}\cdot\nu)\dx+
\frac{\alpha_*C_1^2\, t^2}{2\abs{B}^{1/n}} +
\frac{\alpha'\Vert\varphi-\psi\Vert^2_{L^\infty(\R^n)}}{2\abs{B}^{1/n}},\notag
\end{align}
for some constant $\alpha'>0$ assuming that, say, e.g., 
$\abs{\partial B}\le 1$.
We will eventually choose $\abs{B}< \delta<1$  sufficiently small
(maybe dependent on $\Vert \sigma\Vert_{L^{\infty}({\mathcal R})}$)
and $\Lambda$ sufficiently large 
(independent of $\Vert \sigma\Vert_{L^{\infty}({\mathcal R})}$)
so that
\begin{equation}\label{E6}
C_1^2\,\abs{B}\, t^2\Vert \sigma\Vert_{L^{\infty}({\mathcal R})}
+\frac{\alpha_*C_1^2 t^2}{2\abs{B}^{1/n}}
\le\frac{\Lambda t^2}{2\abs{B}^{1/n}}
\Vert \sigma\Vert_{L^{\infty}({\mathcal R})}.
\end{equation}
Furthermore,
$$
\frac{K_9 \abs{B}}{2} \norm{\sigma}_{L^{\infty}({\mathcal R})}
+\int_{\partial B}\alpha z_\theta (\nabla z_{\theta}\cdot\nu)\dx
\le \frac{\Lambda}{2}\abs{B}\norm{\sigma}_{L^{\infty}({\mathcal R})}.
$$
There exists $\omega>0$ (independent of 
$\Vert \sigma \Vert_{L^{\infty}({\mathcal R})}$) such that
\be\label{E7}
\frac{\alpha_*}{4\Lambda_0\abs{B}^{2/n}}-
K_8-\frac{K_7^2}{\alpha_*}-
\frac{\Vert\nabla\alpha\Vert^2_{L^\infty}}{2\alpha_*}
\ge \frac{\omega}{2\abs{B}^{2/n}}.
\ee
Substituting \eqref{E6} and \eqref{E7} in \eqref{E5}, we have
\begin{align}
\frac{d}{dt}\int_{B}z_\theta^2(t,x)\dx&\le
-\frac{\omega}{\abs{B}^{2/n}}\int_{B}z_\theta^2(t,x)\dx+ \Lambda\abs{B}\,\norm{\sigma}_{L^{\infty}({\mathcal R})}\notag\\
&\quad
+\frac{\Lambda  t^2}{\abs{B}^{1/n}}\norm{\sigma}_{L^{\infty}({\mathcal R})}
+\frac{\alpha'}{\abs{B}^{1/n}}\norm{\varphi-\psi}^2_{L^\infty(\R^n)}.\label{E71}
\end{align}
By the Gronwall inequality and 
\eqref{teta'0}, we have
\begin{align}
\int_{B}z_\theta^2(t,x)\dx&\le
\exp\big(-\frac{\omega t}{\abs{B}^{2/n}}\big)\int_{B}z_\theta^2(0,x)\dx\notag\\
&\quad+\Lambda \exp\big(-\frac{\omega t}{\abs{B}^{2/n}}\big)\int_{0}^{t}
\exp\big(\frac{\omega \tau}{\abs{B}^{2/n}}\big)\abs{B}\,
\Vert \sigma\Vert_{L^{\infty}({\mathcal R})}d\tau\notag\\
&\quad+\exp\big(-\frac{\omega t}{\abs{B}^{2/n}}\big)\int_{0}^{t}
\exp\big(\frac{\omega \tau}{\abs{B}^{2/n}}\big)
\frac{\Lambda  \tau^2}{\abs{B}^{1/n}}
\Vert \sigma\Vert_{L^{\infty}({\mathcal R})}d\tau\notag\\
&\quad+ \exp\big(-\frac{\omega t}{\abs{B}^{2/n}}\big)\int_{0}^{t}
\exp\big(\frac{\omega \tau}{\abs{B}^{2/n}}\big)
\frac{\alpha'}{\abs{B}^{1/n}}
\norm{\varphi-\psi}^2_{L^\infty(\R^n)}d\tau\notag\\
&\le \exp\big(-\frac{\omega t}{\abs{B}^{2/n}}\big)\int_{B}
 \big(\varphi(x)-\psi(x)\big)^2\dx\label{E8}\\
&\quad+\Lambda\frac{\abs{B}^{1+2/n}}{\omega}\Vert 
\sigma\Vert_{L^{\infty}({\mathcal R})}
\left(1-\exp\big(-\frac{\omega t}{\abs{B}^{2/n}}\big)\right)\notag\\
&\quad+
\frac{\Lambda  t^2}{\omega}
\Vert \sigma\Vert_{L^{\infty}({\mathcal R})}\abs{B}^{1/n}
\left(1-\exp\big(-\frac{\omega t}{\abs{B}^{2/n}}\big)\right)\notag\\
&\quad+\frac{\alpha'}{\omega}\norm{\varphi-\psi}^2_{L^\infty(\R^n)}
\abs{B}^{1/n}
\left(1-\exp\big(-\frac{\omega t}{\abs{B}^{2/n}}\big)\right).\notag
\end{align}
Observe that,
\begin{align*}
\frac{\abs{B}^{2/n}}{\omega}
\left(1-\exp\big(-\frac{\omega t}{\abs{B}^{2/n}}\big)\right)&\le t,\quad
1-\exp\big(-\frac{\omega t}{\abs{B}^{2/n}}\big)\le 1,\\
\exp\big(-\frac{\omega t}{\abs{B}^{2/n}}\big)&\le 1,\qquad t\ge 0,
\end{align*}
and, by \eqref{ass1}, \eqref{ass2} and Remark \ref{bounds'},
\begin{align*}
\Vert \sigma\Vert_{L^{\infty}({\mathcal R})}&\le K_0\big(
\Vert a-b\Vert_{L^{\infty}({\mathcal R})}+
\Vert \nabla_x\cdot f-\nabla_x\cdot g\Vert_{L^{\infty}({\mathcal R})}\\
&\quad+\Vert f_u-g_u\Vert_{L^{\infty}({\mathcal R})}+ \Vert h-k\Vert_{L^{\infty}({\mathcal R})}
\big),\qquad 0\le t\le T,
\end{align*}
for some positive constant $K_0$, then, from \eqref{E8} and since $\abs{B}<1$,
\begin{align}
\int_{B}z_\theta^2(t,x)\dx&\le\int_{B}
 \big(\varphi(x)-\psi(x)\big)^2\dx\notag\\
&\quad+
K_0\Lambda\left(1+\frac{1}{\omega}\right)
\abs{B}^{1/n}\,t^2\big(
\Vert a-b\Vert_{L^{\infty}({\mathcal R})}+
\Vert \nabla_x\cdot f-\nabla_x\cdot g\Vert_{L^{\infty}({\mathcal R})}\notag\\
&\quad+ \Vert f_u-g_u\Vert_{L^{\infty}({\mathcal R})}
+ \Vert h-k\Vert_{L^{\infty}({\mathcal R})}\big)
+\frac{\alpha' \abs{B}^{1/n}}{\omega}\norm{\varphi-\psi}^2_{L^\infty(\R^n)}\notag\\
&\le \left(\abs{B}+\frac{\alpha' }{\omega}\abs{B}^{1/n}\right)
\norm{\varphi-\psi}^2_{L^\infty(\R^n)}\notag\\
&\quad+K_0\Lambda\left(1+\frac{1}{\omega}\right)
\abs{B}^{1/n}\,t^2\big(
\Vert a-b\Vert_{L^{\infty}({\mathcal R})}+
\Vert \nabla_x\cdot f-\nabla_x\cdot g\Vert_{L^{\infty}({\mathcal R})}\notag\\
&\qquad+ \Vert f_u-g_u\Vert_{L^{\infty}({\mathcal R})}
+ \Vert h-k\Vert_{L^{\infty}({\mathcal R})}\big).\label{E9}
\end{align}
Let now $\tilde E\supset E$ be a connected set such that interior of
$\tilde E$ contains the closure of $E$, $\text{dist}(\partial\tilde
E,\partial E)>0$, and $\vert\tilde E\vert=2\abs{E}$.  Since the
closure of $E$ is compact, we can cover it with finitely many balls
$B_1,\dots, B_m\subset\R^n$, that is, $E\subset\cup_j B_j$. We may
choose the balls such that $\cup_j B_j$ is contained in the interior
of $\tilde E$, and thus
$$
\abs{\cup_j B_j}\le \sum_{j=1}^k \abs{B_j}\le \vert\tilde E\vert=2\abs{E}.
$$
We assume that both $\abs{\partial B_j}\le 1$ and  $\abs{B_j}\le \delta<1$.  Thus the result \eqref{E9} holds and we may sum the inequality over all balls $B_1,\dots, B_m\subset\R^n$, which yields
\begin{align}
\int_{E}z_\theta^2(t,x)\dx&\le
C\left(\abs{E}+\abs{E}^{1/n}\right)
\norm{\varphi-\psi}^2_{L^\infty(\R^n)}\notag\\
&\quad+C\abs{E}^{1/n}\,t^2\big(
\Vert a-b\Vert_{L^{\infty}({\mathcal R})}+
\Vert \nabla_x\cdot f-\nabla_x\cdot g\Vert_{L^{\infty}({\mathcal R})}\notag\\
&\qquad+ \Vert f_u-g_u\Vert_{L^{\infty}({\mathcal R})}
+ \Vert h-k\Vert_{L^{\infty}({\mathcal R})}\big)\label{E9*}
\end{align}
which proves \eqref{En21}.
\end{proof}

This proves the following result.
\begin{theorem}\label{thmain_p=2} 
Fix $T>0$. Let $u=u(t,x)$ and $v=v(t,x)$ be the classical solution of \eqref{1} and
\eqref{2}, respectively, with $a=a(t,x,y,q)$ and
$b=b(t,x,y,q)$ satisfying $({\mathcal H}_1)$,
$f=f(t,x,y)$, $g=g(t,x,y)$, $h=h(t,x,y,q)$,
and $k=k(t,x,y,q)$ satisfying $({\mathcal H}_2)$, and
$\varphi$ and $\psi$ satisfying $({\mathcal H}_3)$.
 Then there exists a positive constant $C$ depending only on $T,\> n,\> a_*,\>a^*,\>b_*,\>b^*,\> k_1$, $k_2$, and $k_3$
such that
\begin{align}
\Vert u(t,\cdot)-v(t,\cdot)\Vert_{L^2(E)}&\le
C\big(\abs{E}^{1/2n}+\abs{E}^{1/2} \big)
\Vert \varphi-\psi\Vert_{L^\infty(\R^n)}\notag\\
&\quad+ C\,t\,\Big(\Vert a-b\Vert_{L^{\infty}({\mathcal R_0})}+
\Vert \nabla_x\cdot f-\nabla_x\cdot g\Vert_{L^{\infty}({\mathcal R})}
\label{stab_p=2}\\
&\qquad\quad+ \Vert f_u-g_u\Vert_{L^{\infty}({\mathcal R})}
+ \Vert h-k\Vert_{L^{\infty}({\mathcal R_0})}
\Big)^{1/2} \abs{E}^{1/2n},\notag
\end{align}
for all 
$0\le t\le T$ with
${\mathcal R}\doteq [0,T]\times E\times[-K_1,K_1]$, 
${\mathcal R_0}\doteq [0,T]\times E\times[-K_1,K_1]\times[-K_2,K_2]$ 
where $E\subset \R^n$ is bounded connected set with 
Lipschitz boundary.
\end{theorem}
\begin{proof}
Direct consequence of \eqref{curv1},
\eqref{curv2} and Lemmas \ref{En1},  \ref{En2},  and  \ref{En3}.
\end{proof}

\section{Estimates in $L^p(E)$}\label{section3}

We want to extend the estimate of Theorem \ref{thmain_p=2} to general $p$.
\begin{lemma}[{\bf Case $\mathbf {1\le p<2}$}] \label{En2}  Assume  $({\mathcal H}_1)$, $({\mathcal H}_2)$, and $({\mathcal H}_3)$.
There exists a positive constant $C_3$ depending only on 
$T$, $n,\> a_*,\>a^*,\>b_*,\>b^*,\> k_1,\>k_2$ and $k_3$
such that
\begin{align}
\Vert z_\theta(t,\cdot)\Vert_{L^p(E)}&\le
C_3\left(\abs{E}^{(2-p)/(2p)+1/2n}+\abs{E}^{1/p}\right)
\Vert \varphi-\psi\Vert_{L^\infty(\R^n)}\notag\\
&\quad+C_3\, t\,\Big(\Vert a-b\Vert_{L^{\infty}({\mathcal R_0})}+
\Vert \nabla_x\cdot f-\nabla_x\cdot g\Vert_{L^{\infty}({\mathcal R})}
\label{En23}\\
&\quad+ \Vert f_u-g_u\Vert_{L^{\infty}({\mathcal R})}
+ \Vert h-k\Vert_{L^{\infty}({\mathcal R_0})}
\Big)^{1/2} \abs{E}^{(2-p)/(2p)+1/2n},\notag
\end{align}
for each $0\le t\le T$, $E\subset \R^n$ bounded connected set with 
Lipschitz boundary, $0\le \theta \le 1$ and $1\le p<2$.
\end{lemma}
\begin{proof}  By the H\"older inequality,
\begin{align}
\Vert z_\theta(t,\cdot)\Vert_{L^p(E)}^p&=\int_E
z^p_\theta (t,x)\dx\notag\\
&\le \abs{E}^{1/q'}\Vert z_\theta^p(t,\cdot)\Vert_{L^q(E)}\label{E10}\\
&=\abs{E}^{1/q'}\left(\int_E
z^{pq}_\theta (t,x)\dx\right)^{1/q},\notag
\end{align}
with
$$
q\doteq \frac{2}{p},\qquad q'\doteq\frac{2}{2-p}.
$$
So, by \eqref{E10},
\begin{align*}
\Vert z_\theta(t,\cdot)\Vert_{L^p(E)}^p&\le
\abs{E}^{(2-p)/2}\left(\int_E
z^2_\theta (t,x)\dx\right)^{p/2}\\
&\le\abs{E}^{(2-p)/2}\Vert z_\theta(t,\cdot)\Vert_{L^2(E)}^p,
\end{align*}
then, by Lemma \ref{En1},
\begin{align*}
\Vert z_\theta(t,\cdot)\Vert_{L^p(E)}&\le
C_2\left(\abs{E}^{(2-p)/(2p)+1/2n}+\abs{E}^{1/p}\right)\Vert \varphi-\psi\Vert_{L^\infty(\R^n)}\\
&\quad+C_3 t\Big(\Vert a-b\Vert_{L^{\infty}({\mathcal R})}+
\Vert \nabla f-\nabla g\Vert_{L^{\infty}({\mathcal R})}\\
&\quad+ \Vert f_u-g_u\Vert_{L^{\infty}({\mathcal R})}
+ \Vert h-k\Vert_{L^{\infty}({\mathcal R})}
\Big)^{1/2} \abs{E}^{(2-p)/(2p)+1/2n},
\end{align*}
This concludes the proof.
\end{proof}
\begin{lemma}[{\bf Case $\mathbf{p>2}$}] \label{En3}Assume  
$({\mathcal H}_1)$, $({\mathcal H}_2)$, and $({\mathcal H}_3)$.
There exists a positive constant $C_4$ depending only on $T$,
$n,\> a_*,\>a^*,\>b_*,\>b^*,\> k_1,\> k_2$ and $k_3$
such that
\begin{align}
\Vert z_\theta(t,\cdot)\Vert_{L^p(E)}&\le
C_4(1+ t^{(p-2)/p})\big(\abs{E}^{1/np}+\abs{E}^{1/p}\big)\Vert \varphi-\psi\Vert_{L^\infty(\R^n)}^{2/p}\notag\\
&\quad+C_4 (t+t^{2/p})
\Big(\Vert a-b\Vert_{L^{\infty}({\mathcal R})}+
\Vert \nabla_x\cdot f-\nabla_x\cdot g\Vert_{L^{\infty}({\mathcal R})}
\label{En4}\\
&\quad+ \Vert f_u-g_u\Vert_{L^{\infty}({\mathcal R})}
+ \Vert h-k\Vert_{L^{\infty}({\mathcal R})}
\Big)^{1/p} \abs{E}^{1/(np)},\notag
\end{align}
for each $0\le t\le T$, $E\subset \R^n$ bounded connected set 
with Lipschitz boundary, $0\le \theta \le 1$ and $ 2<p<\infty$.
\end{lemma}
\begin{proof} Observe that
\begin{align*}
\Vert z_\theta(t,\cdot)\Vert_{L^p(E)}^p &=\int_E
z^p_\theta (t,x)\dx\\
&\le\Vert z_\theta(t,\cdot)\Vert_{L^\infty(\R^n)}^{p-2}\int_E
z^2_\theta (t,x)\dx\\
&=\Vert z_\theta(t,\cdot)\Vert_{L^\infty(\R^n)}^{p-2}
\Vert z_\theta(t,\cdot)\Vert_{L^2(E)}^{2}.
\end{align*}
Since $2/p,\> (p-2)/p<1$, by Lemmas \ref{teta1} and \ref{En1},
we have
\begin{align*}
\Vert z_\theta(t,\cdot)\Vert_{L^p(E)}&\le
\Vert z_\theta(t,\cdot)\Vert_{L^\infty(\R^n)}^{(p-2)/p}\,
\Vert z_\theta(t,\cdot)\Vert_{L^2(E)}^{2/p}\\
&\le\Big(C_1^{(p-2)/p}\,t^{(p-2)/p}
+\Vert \varphi-\psi \Vert_{L^{\infty}(\R^n)}^{(p-2)/p}\Big)\\
&\quad \times\Big[C_2^{2/p}\big(\abs{E}^{1/np}+\abs{E}^{1/p} \big)
\Vert \varphi-\psi \Vert_{L^{\infty}(\R^n)}^{2/p}\\
&\quad+C_2^{2/p}\,t^{2/p}\abs{E}^{1/np}
\Big(\Vert a-b\Vert_{L^{\infty}({\mathcal R})}+
\Vert \nabla_x\cdot f-\nabla_x\cdot g\Vert_{L^{\infty}({\mathcal R})}\\
&\quad+ \Vert f_u-g_u\Vert_{L^{\infty}({\mathcal R})}
+ \Vert h-k\Vert_{L^{\infty}({\mathcal R})}\Big)^{1/p}\Big] \\
&\le \Big(C_1^{(p-2)/p}\,t^{(p-2)/p}+k_4\Big) \\
&\quad
\times\Big[C_2^{2/p}\big(\abs{E}^{1/np}+\abs{E}^{1/p}\big)
\Vert \varphi-\psi \Vert_{L^{\infty}(\R^n)}^{2/p}\\
&\quad
+C_2^{2/p}\,t^{2/p}\abs{E}^{1/np}
\Big(\Vert a-b\Vert_{L^{\infty}({\mathcal R})}+
\Vert \nabla_x\cdot f-\nabla_x\cdot g\Vert_{L^{\infty}({\mathcal R})}\\
&\quad+ \Vert f_u-g_u\Vert_{L^{\infty}({\mathcal R})}
+ \Vert h-k\Vert_{L^{\infty}({\mathcal R})}\Big)^{1/p}\Big] \\
&=\Big(C_1^{(p-2)/p}\,t^{(p-2)/p}+k_4\Big)C_2^{2/p}
\big(\abs{E}^{1/np}+\abs{E}^{1/p} \big)\\
&\quad\times\Vert \varphi-\psi \Vert_{L^{\infty}(\R^n)}^{2/p}
 +C_2^{2/p}\Big(C_1^{(p-2)/p}\,t+k_4t^{2/p}\Big)\abs{E}^{1/np}\\
&\quad \times\Big(\Vert a-b\Vert_{L^{\infty}({\mathcal R})}+
\Vert \nabla_x\cdot f-\nabla_x\cdot g\Vert_{L^{\infty}({\mathcal R})}\\
&\quad+ \Vert f_u-g_u\Vert_{L^{\infty}({\mathcal R})}
+ \Vert h-k\Vert_{L^{\infty}({\mathcal R})}\Big)^{1/p},
\end{align*}
where $k_4$ is a positive constant such that
$$
(2k_3)^{(p-2)/p}\le k_4,\qquad 2<p<\infty.
$$
Since the maps
$$
2<p<\infty\longmapsto C_1^{(p-2)/p},\> C_2^{2/p}
$$
are bounded the proof is done. 
\end{proof}

The following theorem summarizes the result in Theorem \ref{thmain_p=2} with 
the extension to general $p$.
\begin{theorem}\label{thmain} 
Fix $T>0$. Let $u=u(t,x)$ and $v=v(t,x)$ be the classical solution of
\eqref{1} and \eqref{2}, respectively, with $a=a(t,x,y,q)$ and
$b=b(t,x,y,q)$ satisfying $({\mathcal H}_1)$, $f=f(t,x,y)$,
$g=g(t,x,y)$, $h=h(t,x,y,q)$, and $k=k(t,x,y,q)$ satisfying $({\mathcal
H}_2)$, and $\varphi$ and $\psi$ satisfying $({\mathcal H}_3)$.  Then
there exists a positive constant $C$ depending only on $T,\> n,\>
a_*,\>a^*,\>b_*,\>b^*,\> k_1$, $k_2$, and $k_3$ such that
\begin{align}
\Vert u(t,\cdot)-v(t,\cdot)\Vert_{L^p(E)}&\le
A_E(t)\Vert \varphi-\psi\Vert^{2\rho_p}_{L^\infty(\R^n)}\notag\\
&\quad+ B(t)\Big(\Vert a-b\Vert_{L^{\infty}({\mathcal R})}+
\Vert \nabla_x\cdot f-\nabla_x\cdot g\Vert_{L^{\infty}({\mathcal R})}
\label{stab}\\
&\qquad+ \Vert f_u-g_u\Vert_{L^{\infty}({\mathcal R})}
+ \Vert h-k\Vert_{L^{\infty}({\mathcal R})}
\Big)^{\rho_p} \abs{E}^{\eta_p},\notag
\end{align}
with
${\mathcal R}\doteq [0,T]\times E\times[-K_1,K_1]$, 
${\mathcal R_0}\doteq [0,T]\times E\times[-K_1,K_1]\times[-K_2,K_2]$. 
Here
\begin{align*}
\rho_p&\doteq\begin{cases}
\frac{1}{2},  &\text{if $1\le p\le 2$},  \\
\frac{1}{p},  &\text{if $2< p<\infty$},
\end{cases}
\qquad
\eta_p\doteq\begin{cases}
\frac{2-p}{2p}+\frac{1}{2n},  &\text{if $1\le p\le 2$},  \\
\frac{1}{np},  &\text{if $2< p<\infty$},
\end{cases}\\
 A_E(t)&\doteq C\begin{cases}
(\abs{E}^{(2-p)/2p+1/2n}+\abs{E}^{1/p}),  &\text{if $1\le p\le 2$},  \\
(1+t^{(p-2)/p})(\abs{E}^{1/np}+\abs{E}^{1/p}),  &\text{if $2< p<\infty$},
\end{cases}\\
B(t)&\doteq C\begin{cases}
t,  &\text{if $1\le p\le 2$},  \\
(t+t^{2/p}),  &\text{if $2< p<\infty$},
\end{cases}
\end{align*}
for all $0\le t\le T$,  where $E\subset \R^n$ is bounded connected set with Lipschitz boundary and $1\le p<\infty$.
\end{theorem}
\begin{proof}  
Direct consequence of \eqref{curv1},
\eqref{curv2} and Lemmas \ref{En1},  \ref{En2}, \ref{En3}.
\end{proof}

\bigskip
\noindent {\bf Acknowledgments.}  The authors
would like to thank Prof.\ Alberto Bressan many useful discussions.

%
%

\end{document}